\documentclass[12pt]{amsart}

\usepackage{amsmath}
\usepackage{amssymb}
\usepackage{amsthm}

\usepackage[dvips]{graphics}
\usepackage{graphicx}
\usepackage{epsfig}

\usepackage{ifthen}

\newcommand{\cl}{\operatorname{cl}}

\newcommand{\cof}{\operatorname{cof}}

\newcommand{\Rr}{{\mathbb{R}}}

\newcommand{\Zz}{{\mathbb{Z}}}

\newcommand{\Tt}{{\mathbb{T}}}

\newcommand{\Hh}{{\overline{H}}}

\newcommand{\Cg}{{\mathcal{C}}}

\newcommand{\Ii}{{\mathcal{I}}}

\newcommand{\Mm}{{\mathcal{M}}}
\newcommand{\Mn}{{\mathcal{N}}}

\newcommand{\cqd}{\hfill $\blacksquare$}
\newcommand{\cqda}{\hfill $\rule{1.5mm}{1.5mm}$}

\newcommand{\Pf}{{\noindent \sc Proof. \ }}

\newcommand{\cws}{{\overset{*}{\rightharpoonup}}}

\newtheorem{teo}{Theorem}
\newtheorem{lem}{Lemma}
\newtheorem{pro}{Proposition}
\newtheorem{cor}{Corollary}

\theoremstyle{definition}

\begin{document}
\title{A Stochastic Analog of Aubry-Mather theory}
\author{Diogo Aguiar Gomes} 

\begin{abstract}
In this paper we discuss a stochastic analog of
Aubry-Mather theory in which a deterministic control problem
is replaced by a controlled diffusion. We prove the
existence of a minimizing measure (Mather measure) and discuss
its main properties using viscosity solutions
of Hamilton-Jacobi equations. Then
we prove regularity estimates
on viscosity solutions of Hamilton-Jacobi equation
using the Mather measure.
Finally we apply these results to
prove asymptotic estimates on the 
trajectories of controlled diffusions and
study the convergence of Mather measures
as the rate of diffusion vanishes.
\end{abstract}

\maketitle
\tableofcontents

\section{Introduction}

The objective of this paper is to understand a stochastic analog
of Aubry-Mather theory. The original problem \cite{Mather} consists
in determining a probability measure $\mu$ in $\Tt^n\times \Rr^n$
($\Tt^n$ is the $n$-dimensional torus)
that
minimizes the average action
\begin{equation}
\label{one}
\int L(x,v)d\mu
\end{equation}
for a given Lagrangian $L$ with the constraint that $\mu$ is
invariant under the flow generated by
the Euler-Lagrange equations
associated with $L$. This problem is equivalent
\cite{MatherPC}
to the relaxed problem of
minimizing (\ref{one}) with the constraint
$$
\int vD_x\phi d\mu=0
$$
for any $\phi(x)$. In the case of controlled diffusions,
we replace this constraint by
$$
\int A^v\phi d\mu=0,
$$
for all $\phi$ smooth and periodic,
in which $A^v$ is the infinitesimal generator of the controlled
diffusion.

We proceed as follows: in section \ref{sec2}
we construct a relaxed
minimization problem on a space of measures.
Then, in section \ref{dois},
we identify its dual by means of
Fenchel-Rockafellar
duality theorem \cite{Rock}. This
dual problem turns out to involve
Hamilton-Jacobi equations which are studied in section \ref{ApA}.
We prove equivalence between the
strong and weak problems (section \ref{seceq}), and
characterize the minimizing measures
using viscosity solutions of Hamilton-Jacobi equations (section \ref{sec6}).
The we discuss several
applications:  regularity
of Hamilton-Jacobi equations (section \ref{sec7}),
logarithmic transform, connection with 
eigenvalue problems (section \ref{sec8})
asymptotics for controlled diffusions (section \ref{sec9})
and convergence of the stochastic Mather measure
as the diffusion coefficient vanishes (section \ref{sec10}).

The original Mather problem, as well as its stochastic version
are convex linear
programming problem over a space of Radon measures.
Related control problems have been studied by duality
\cite{Vinter3},
\cite{Vinter2}, \cite{Vinter1}, 
\cite{MR90i:49016}, \cite{MR89k:93223} and \cite{MR91a:49003},
in which
Fenchel-Rockafellar duality theorem \cite{Rock} is used to analyze optimal
control problems. In this paper we apply similar techniques to understand
Aubry-Mather theory and its
stochastic analogs.

Several authors have studied the relation between viscosity solutions
of Hamilton-Jacobi equations and Mather measures \cite{FATH1}, \cite{FATH2},
\cite{FATH3}, \cite{FATH4}, \cite{E}, \cite{EGom1}, \cite{Gomes3} and
\cite{Gomes2}.
The results by
A. Fathi
\cite{FATH1}, \cite{FATH2}, \cite{FATH3},
\cite{FATH4},
and W. E \cite{E} make clear
the connection between
viscosity solutions and Hamiltonian dynamics.
The main idea is that if
$u(x,P)$ is a viscosity solution of 
\begin{equation}
\label{two}
H(P+D_xu,x)=\Hh
\end{equation}
(here $H$ is the Legendre transform of $L$)
then
there
exists an invariant set $\Ii$ contained on the graph
$$
\{(x,P+D_xu(x,P))\}.
$$
Furthermore, $\Ii$ is a subset of a Lipschitz graph, i.e.,
$D_xu(x,P)$ is a Lipschitz
function on $\pi(\Ii)$, where $\pi(x,p)=x$.
If $\Hh$ is differentiable at $P$, then any solution 
$(x(t),p(t))$ of (\ref{one}) with initial
conditions on $\Ii$ satisfies
\begin{equation}
\label{12}
\lim_{t\rightarrow\infty}\frac{|x(t)+D_P\Hh t|}{t}=0.
\end{equation}
In \cite{Gomes3} and \cite{Gomes2} this problem is studied
with detail and
more precise asymptotic results are presented.
We also prove 
regularity results for viscosity solution of (\ref{two}) - in
particular uniform continuity in $P$. 
In \cite{EGom1} Mather measures are used
to prove regularity for
solutions of Hamilton-Jacobi equations. The main results
are $L^2$ type estimates in the difference quotients of $D_xu$.
The objective of this paper is to generalize these results to
the stochastic case. 

\section{Stochastic Mather measures}
\label{sec2}

In this section we define a stochastic analog of the Mather's
minimal measure problem \cite{Mather3}, \cite{Mather},
\cite{Mane2}, \cite{Mane5}.
To do so we consider an ergodic diffusion control 
problem and
study an associated relaxed
minimization problem on a space of measures.
In the next section,
we identify its dual by means of
Fenchel-Rockafellar
duality theorem and show that the
dual problem is, in some sense, a
Hamilton-Jacobi equation.
%Then, in the
%remaining sections, we prove the equivalence between the
%strong and weak problems, characterize the minimizing measures
%using viscosity solutions of Hamilton-Jacobi equations, and discuss several
%applications - convergence as the diffusion coefficient vanishes, regularity
%of Hamilton-Jacobi equations, and asymptotics for controlled diffusions.

Consider a controlled Markov diffusion \cite{FS} in $\Rr^n$
\begin{equation}
\label{diff}
dx=v dt+\sigma dw,
\end{equation}
where $v$ is a progressively measurable control, $w$ a $n$-dimensional
Brownian motion and $\sigma\geq 0$
the diffusion rate ($\sigma=0$ corresponds to the standard Aubry-Mather
theory).
The control objective
is to minimize the long-time running cost
$$
\lim_{T\rightarrow+\infty}\frac{1}{T}E \int_0^T L(x,v)dt,
$$
over all admissible control processes $v$, this is 
called the ergodic control problem (here $E$ denotes
the expected value with respect to the underlying probability measure).
We assume that
the function $L(x,v)$ is smooth in both variables, $\Zz^n$ periodic in $x$,
coercive
and strictly convex in $v$. Furthermore, since adding a constant to $L$ does
not change the nature of the problem, we also assume $L\geq 0$.

Let $\Omega=\Tt^n\times \Rr^n$, where $\Tt^n$ is the $n$-dimensional torus,
identified, when convenient, with $[0,1]^n$ or $\Rr^n$ with a
periodic structure (in geometric terms, $\Rr^n$ is the universal covering
of $\Tt^n$). A pair
$(x,v)=z$ represents
a generic point $z\in\Omega$, with $x\in \Tt^n$ and $v\in\Rr^n$.
Choose a function
$\gamma\equiv\gamma(|v|):\Omega\rightarrow [1,+\infty)$
satisfying
$$
\lim_{|v|\rightarrow +\infty}\frac{L(x,v)}{\gamma(v)}=+\infty\qquad
\lim_{|v|\rightarrow +\infty}\frac{|v|}{\gamma(v)}=0.
$$
Let $\Mm$ be the set of weighted Radon measures on $\Omega$, i.e.,
$$
\Mm=\{\text{signed measures on }\Omega \text{ with }
\int_\Omega \gamma d|\mu|<\infty\}.
$$
Note that $\Mm$ is the dual of
the set $C_\gamma^0(\Omega)$ of continuous
functions $\phi$ with
$$
\|\phi\|_\gamma=\sup_{\Omega} |\frac{\phi}{\gamma}|<\infty, \qquad
\lim_{|z|\rightarrow \infty}\frac{\phi(z)}{\gamma(v)}\rightarrow 0.
$$
For each bounded control strategy $v$ consider the measure $\mu_T$
defined by
$$
\int \phi d\mu_T=\frac{1}{T}E\int_0^T \phi(x(t),v(t))dt.
$$
As $T\rightarrow +\infty$ we may extract a weakly convergent
subsequence $\mu_T\rightharpoonup \mu(v)$. Let
$$
\Mm_0=\cl \{\mu(v):v\ \text{bounded control strategy}\}.
$$
Define
$$
\Mm_1=\{\mu\in\Mm :\int_\Omega d\mu = 1, \mu\geq 0\}.
$$
The stochastic analog of Mather's problem consist in determine
a measure $\mu$ that minimizes
\begin{equation}
\label{n1}
\inf_{\mu \in \Mm_0\cap\Mm_1} \int_\Omega Ld\mu.
\end{equation}
For our purposes, however,
it is convenient to consider a relaxed problem by replacing $\Mm_0$ by
a slightly larger set $\Mn_0$ that we define next.

The infinitesimal generator corresponding to the controlled diffusion
(\ref{diff})
is
$$
A^v\phi=\frac{\sigma^2}{2}\Delta \phi+v\cdot \nabla\phi.
$$
\begin{pro}
Any measure $\mu(v)$ in $\Mm_0$ satisfies
\begin{equation}
\label{admit}
\int A^v\varphi d\mu=0,
\end{equation}
for all $\varphi=\varphi(x)$,
$\varphi$ periodic and $C^2$ (or $C^1$ if $\sigma=0$).
\end{pro}
\Pf
Consider the measure defined by
$$
\int_\Omega \phi d\mu_T=\frac{1}{T}E\int_0^T \phi(x(t),v(t))dt.
$$
Assume $\mu_T\cws \mu$. We claim that
$$
\int_\Omega A^v\varphi d\mu=0,
$$
for $\varphi(x)$ $C^2$ and periodic function of $x$ only.
To see this
recall
Dynkin's formula:
$$
\varphi(x(T))-\varphi(x(0))=E\int_0^T A^{v(t)}\varphi (x(t))dt
$$
for any $x(t)$ and $v(t)$ that solve (\ref{diff}).
In the case $\sigma=0$ this is just the fundamental theorem
of calculus.
Dividing by $T$ and
letting $T\rightarrow \infty$ we obtain (\ref{admit}).
\cqd

Let $\Mn_0$ be the
closure of the set of all measures that satisfy (\ref{admit}):
$$
\Mn_0=\cl\{\mu\in \Mm:\int_\Omega A^v\varphi d\mu=0,\,
\forall \varphi(x)\in C^2(\Tt^n)\}.
$$
In the case $\sigma=0$, the set $\Mn_0$ is the ``measure theoretic''
analog of the set of
closed curves on $\Tt^n$. Indeed, if $\theta:[0,1]\rightarrow \Tt^n$ is
a piecewise smooth closed curve the we can define a measure $\mu_\theta$ by
$$
\int_\Omega fd\mu_\theta=\int_0^1 f(\theta(t),\dot \theta(t))dt.
$$
Clearly $\mu_\theta$ is in $\Mn_0$, and since $\Mn_0$ is
a linear space, it contains all linear combinations of
measures of this form.

The additional problem that we will consider is
\begin{equation}
\label{n2}
\inf_{\mu\in \Mn_0\cap\Mm_1} \int Ld\mu.
\end{equation}
We will prove later on that
\begin{equation}
\label{ident}
\inf_{\mu\in \Mn_0\cap\Mm_1} \int Ld\mu=\inf_{\mu\in \Mm_0\cap\Mm_1}
\int Ld\mu.
\end{equation}
This identity is a consequence that
$\Mn_0$
is the weak-$*$ closure of
$\Mm_0$. However, the proof of this depends on
(\ref{ident}) holding for a sufficiently large class of $L$ (see
\cite{MR90i:49016}, \cite{MR89k:93223}, and \cite{MR91a:49003} for
related proofs). Therefore we will prove (\ref{ident}) directly.

The last
issue we discuss
in this section is the existence of a measure that minimizes:
$$
\inf_{\mu\in\Mm_0\cap\Mm_1}\int Ld\mu.
$$
This measure is the stochastic analog of the Aubry-Mather measure.
A similar proof also shows that there exists a minimizing measure
in $\Mn_0\cap \Mm_1$.
In the next section we prove that
$$
\inf_{\mu\in\Mm_0\cap\Mm_1}\int Ld\mu=\inf_{\mu\in\Mn_0\cap\Mm_1}\int Ld\mu.
$$

First we quote a compacity lemma:
\begin{lem}[Ma\~n\'e \cite{Mane5}]
In $\Mm_0\cap\Mm_1$ the set
$$
\int Ld\mu<c
$$
is compact with respect to the weak-$*$ topology in $(C_\gamma^0)'$
\end{lem}

%For a proof of the previous lemma consult \cite{Mane5}.
With the help of this lemma we prove the existence of a minimizing
measure.

\begin{teo}
There exists a measure $\mu\in\Mm_0\cap\Mm_1$ such that:
$$
\int Ld\mu=\inf_{\mu\in\Mm_0\cap\Mm_1}\int Ld\mu.
$$ 
\end{teo}
\Pf Take any minimizing sequence $\mu_n$. Since $\int Ld\mu<c$,
the previous lemma shows that by extracting a subsequence, if necessary,
$\mu_n\cws\mu$. Thus,
for any fixed $k$,
$$
\lim_{n\rightarrow +\infty}
\int \min(L,k)d\mu_n\rightarrow \int \min(L,k)d\mu.
$$ 
Thus
$$
\int \min(L,k)d\mu\leq \inf_{\mu\in\Mm_0\cap\Mm_1}\int Ld\mu,
$$
for all $k$. But then by monotone convergence theorem
$$
\int Ld\mu\leq \inf_{\mu\in\Mm_0\cap\Mm_1}\int Ld\mu,
$$
which proves the theorem.
\cqd

A similar proof yields:
\begin{teo}
There exists a measure $\mu\in\Mn_0\cap\Mm_1$ such that:
$$
\int Ld\mu=\inf_{\mu\in\Mn_0\cap\Mm_1}\int Ld\mu.
$$ 
\end{teo}

\section{Identification of the dual problem}
\label{dois}

In this section we identify the dual problem of
$$
\min_{\mu\in \Mn_0\cap\Mm_1}\int Ld\mu.
$$
The dual problem involves a Hamilton-Jacobi equation. Further analysis
of this equation is carried out in the remaining
sections and
yields
important information about the minimizing measure.

First we review some facts about convex duality.
Let $E$ be a Banach space with dual $E'$. The pairing between
$E$ and $E'$ is
denoted by $(\cdot,\cdot)$.
Suppose $h:E\rightarrow (-\infty,+\infty]$ is
a
convex, lower semicontinuous function.
The Legendre-Fenchel transform $h^*:E'\rightarrow [-\infty,+\infty]$
of $h$ is defined by
$$
h^*(y)=\sup_{x\in E} \left(-(x,y)-h(x)\right),
$$
for $y\in E'$. Similarly, for concave, upper semicontinuous
functions $g:E\rightarrow (-\infty,+\infty]$ let
$$
g^*(y)=\inf_{x\in E} \left(-(x,y)-g(x)\right).
$$
%The Fenchel-Rockafellar duality theorem \cite{Rock} is
%more precisely,
\begin{teo}[Rockafellar \cite{Rock}]
\label{dualteo}
Let $E$ be a locally convex Hausdorff topological vector space over $\Rr$
with dual $E^*$.
Suppose $h:E\rightarrow (-\infty,+\infty]$
is convex and lower semicontinuous, $g:E\rightarrow [-\infty,+\infty)$
is concave and upper semicontinuous. Then
\begin{equation}
\label{RFD}
\sup_x g(x)-f(x)=\inf_y f^*(y)-g^*(y),
\end{equation}
provided
that either
$h$ or $g$ is continuous at some point where both functions
are finite. 
\end{teo}

For $\phi\in C_\gamma^0(\Omega)$ define
$$
h_1(\phi)=\sup_{(x,v)\in \Omega}
(-\phi(x,v)-L(x,v)).
$$
Let $\Cg$ be defined by 
$$
\Cg=\cl\{\phi:\phi=A^v\varphi, \varphi(x)\in  C^2(T^n)\},
$$
here $\cl$ denotes the closure in $C_\gamma^0$ (if $\sigma=0$ we may take
$\varphi(x)\in  C^1(T^n)$).
If $\sigma=0$, we may
think of the elements in $\Cg$ as generalized closed differential
forms; indeed if $\theta:[0,1]\rightarrow T^n$ is
a piecewise smooth closed curve and $\phi\in \Cg$ then
$$
\int \phi d\mu_\theta=0.
$$
Define
$$
h_2(\phi)=
\begin{cases}
0\quad&\text{if } \phi\in\Cg\\
-\infty&\text{otherwise}.
\end{cases}
$$
In this section we prove that
\begin{equation}
\label{63b}
\sup_{\phi\in C_\gamma^0(\Omega)} h_2(\phi)-h_1(\phi)
\end{equation}
is
the dual problem of (\ref{n2}).

First we compute
the Legendre-Fenchel
transforms of $h_1$ and $h_2$
in order to apply theorem \ref{dualteo}
to (\ref{63b}).
\begin{pro}
We have
$$
h_1^*(\mu)=
\begin{cases}
\int Ld\mu \quad &\text{if} \quad \mu\in\Mm_1\\
+\infty  &\text{otherwise,}
\end{cases}
$$
and
$$
h_2^*(\mu)=
\begin{cases}
0 \quad &\text{if} \quad \mu\in\Mn_0\\
-\infty  &\text{otherwise.}
\end{cases}
$$
\end{pro}

\Pf
Recall that
$$
h_1^{*}(\mu)=
\sup_{\phi\in C_\gamma^0(\Omega)} \left(-\int\phi d\mu-h_1(\phi)\right).
$$
We claim that if $\mu$ is non-positive
then $h_1^*(\mu)=\infty$.
\begin{lem}
If $\mu\not \geq 0$ then $h_1^*(\mu)=+\infty$.
\end{lem}
\Pf
If $\mu\not \geq 0$ we can choose a sequence of positive functions
$\phi_n\in C_\gamma^0(\Omega)$ such that
$$
\int -\phi_nd\mu\rightarrow +\infty.
$$
Thus, since $L\geq 0$,
$$
\sup_\Omega -\phi_n-L\leq 0.
$$
Therefore if $\mu\not\geq 0$ then
$h_1^*(\mu)=+\infty$.
\cqda

\begin{lem}
If $\mu\geq 0$ then
$$
h_1^{*}(\mu)\geq \int Ld\mu +\sup_{\psi\in C_\gamma^0(\Omega)}\left(
\int\psi d\mu-\sup \psi
\right).
$$
\end{lem}
\Pf 
Let $L_n$ be a sequence of functions in $C_\gamma^0(\Omega)$ increasing
pointwise to $L$. Any function $\phi$ in $C_\gamma^0(\Omega)$
can be written as
$\phi=-L_n-\psi$, for some $\psi$ also in $C_\gamma^0(\Omega)$.
Thus
\begin{align}
\notag
&\sup_{\phi\in C_\gamma^0(\Omega)}
\left(-\int\phi d\mu-h_1(\phi)\right)=\\\notag
&\qquad=\sup_{\psi\in C_\gamma^0(\Omega)}
\left(
\int L_n d\mu+\int \psi d\mu-\sup (L_n+\psi-L)
\right).
\end{align}
Since $L_n-L\leq 0$,
$$
\sup_\Omega L_n-L\leq 0,
$$
thus
$$
\sup_\Omega(L_n+\psi-L)\leq \sup_\Omega \psi.
$$
Thus
$$
\sup_{\phi\in C_\gamma^0(\Omega)}
\left(-\int\phi d\mu-h_1(\phi)\right)\geq 
\sup_{\psi\in C_\gamma^0(\Omega)}
\left(
\int L_n d\mu+\int \psi d\mu-\sup (\psi)
\right).
$$
By the
monotone convergence theorem $\int L_n d\mu\rightarrow \int Ld\mu$. Therefore
$$
\sup_{\phi\in C_\gamma^0(\Omega)}
\left(-\int\phi d\mu-h_1(\phi)\right)\geq 
\int Ld\mu+\sup_{\psi\in C_\gamma^0(\Omega)}
\left(\int \psi d\mu-\sup (\psi)
\right),
$$
as required.
\cqda

If $\int Ld\mu=+\infty$ then $h_1^*(\mu)=+\infty$.
If $\int d\mu\neq 1$ then
$$
\sup_{\psi\in C_\gamma^0(\Omega)}\left(
\int\psi d\mu-\sup \psi\right)\geq
\sup_{\alpha\in \Rr} \alpha(\int d\mu-1)=+\infty,
$$
by taking $\psi\equiv\alpha$, constant.
Therefore $h_1^*(\mu)=+\infty$.

If $\int d\mu=1$ we have, from the previous lemma
$$
h_1^*(\mu)\geq \int Ld\mu,
$$
by taking $\psi\equiv 0$.

Also, for any $\phi$
$$
\int (-\phi-L)d\mu\leq \sup_\Omega (-\phi-L),
$$
if $\int d\mu=1$.
Hence
$$
\sup_{\phi\in C_\gamma^0(\Omega)}
\left(-\int\phi d\mu-h_1(\phi)\right)\leq \int Ld\mu.
$$
Thus
$$
h_1^*(\mu)=
\begin{cases}
\int Ld\mu \quad &\text{if} \quad \mu\in\Mm_1\\
+\infty  &\text{otherwise.}
\end{cases}
$$

Now we will compute $h_2^*$.
First observe that 
if $\mu\not\in\Mn_0$ then there exists $\hat \phi\in\Cg$ such that
$$
\int \hat \phi d\mu\neq 0.
$$
and so
$$
\inf_{\phi\in\Cg}-\int \phi d\mu\leq \inf_{\alpha\in\Rr}\alpha
\int \hat \phi d\mu=-\infty.
$$
If $\mu\in \Mn_0$ then $\int \phi d\mu=0$, for all $\phi\in \Cg$.
Therefore
$$
h_2^{*}(\mu)=\inf_{\phi\in\Cg}-\int \phi d\mu=
\begin{cases}
0\quad&\text{if }\mu\in\Mn_0\\
-\infty&\text{otherwise}.
\end{cases}
$$
\cqd

The Fenchel-Rockafellar duality theorem states
that
\begin{equation}
\label{n3}
\sup_{\phi\in C_\gamma^0(\Omega)}(h_2(\phi)-h_1(\phi))=\inf_{\mu\in\Mm}
(h_1^{*}(\mu)-h_2^{*}(\mu)),
\end{equation}
provided on the set $h_2>-\infty$, $h_1$ is continuous.
In the next lemma we prove that 
$h_1$ is continuous, and therefore (\ref{n3}) holds.
\begin{lem}
$h_1$ is continuous.
\end{lem}
\Pf
Suppose $\phi_n\rightarrow \phi$ in $C_\gamma^0$.
We must prove that $h_1(\phi_n)\rightarrow h_1(\phi)$.
Observe that $\|\phi_n\|_\gamma$
and $\|\phi\|_\gamma$ are bounded uniformly by some constant $C$.
The growth condition on $L$ implies that there exists $R>0$ such that
$$
\sup_{\Omega} -\hat \phi-L=\sup_{\Tt^n\times B_R}-\hat \phi-L,
$$
for all $\hat \phi$ in $C^\gamma_0(\Omega)$ with $\|\hat\phi\|_\gamma<C$,
in which $B_R=\{v\in \Rr^n: |v|\leq R\}$ is the ball
of radius $R$ centered at the origin.
On $B_R$, $\phi_n\rightarrow\phi$ uniformly and so
$$
\sup_{\Omega} -\phi_n-L\rightarrow 
\sup_{\Omega} -\phi-L.
$$
\cqd

Denote by $H^\star$ the value
$$
H^\star=-\sup_{\phi\in C^\gamma_0(\Omega)}(h_2(\phi)-h_1(\phi))
$$
\begin{teo}
\label{nlm}
We have
$$
H^\star=
\inf \{\lambda: \exists \varphi\in C^1(\Tt^n):
-\frac{\sigma^2}{2}\Delta\varphi
+H(D_x\varphi,x)<\lambda\},
$$
in which
$$
H(p,x)=\sup_v -p\cdot v-L(x,v)
$$
is the Legendre transform of $L$.
\end{teo}
\Pf
Note that
\begin{align}
\notag
H^\star&=
\inf_{\varphi\in C^1(T^n)}\sup_{(x,v)\in\Omega}
-\frac{\sigma^2}{2}\Delta\varphi-
vD_x\varphi-L=\\\notag&=
\inf_{\varphi\in C^1(\Tt^n)}\sup_{x\in \Tt^n}
-\frac{\sigma^2}{2}\Delta\varphi+
H(D_x\varphi,x).
\end{align}
\cqd

%%%%%%%%%%%%%%%%%%%%%%%%%%%%%%%%%%%%%%%%%%%

\section{The Cell Problem}
\label{ApA}

The last theorem in the previous section
suggests that we study
the equation
\begin{equation}
\label{c1}
-\frac{\sigma^2}{2}\Delta u+H(D_xu,x)=\Hh.
\end{equation}
In this section we prove that there exists a unique number
$\Hh$ for which (\ref{c1}) has a periodic viscosity solution.
Using the results from \cite{KrylovE}, we show that such solution
is $C^2$. Then we prove that the solution
is unique (up to additive constants). Finally we 
prove estimates on $\Hh$ and $u$ that do not depend on $\sigma$.

\begin{teo}
\label{solcel}
There exists a unique number $\Hh$ for which the equation
$$
-\frac{\sigma^2}{2}\Delta u+H(D_xu,x)=\Hh
$$
has a periodic viscosity solution. Furthermore the solution is $C^2$ and
unique.
\end{teo}
\Pf
First we address the issue of the existence of a viscosity solution. To do so
consider the infinite horizon discounted cost problem
$$
u^\alpha=\inf E\int_0^\infty e^{-\alpha t} L(x,v)dt
$$
with $dx=vdt+\sigma dw$. Then $u^\alpha$ is
a periodic viscosity solution of \cite{FS}
$$
-\frac{\sigma^2}{2}\Delta u^\alpha +H(D_xu^\alpha,x)+\alpha u^\alpha=0.
$$
Since $u^\alpha$ is periodic, uniformly Lipschitz in $\alpha$
\cite{FS}
there
exists a subsequence $u^\alpha$ and $u$ periodic for which
$$
u^\alpha-\min u^\alpha \rightarrow u.
$$
Since $0\leq u^\alpha\leq \frac{C}{\alpha}$ we have
$\alpha u^\alpha\rightarrow -\Hh$, for some $\Hh$ (extracting a further
subsequence if necessary). Then $u$ is a periodic
viscosity solution of
$$
-\frac{\sigma^2}{2}\Delta u +H(D_xu,x)=\Hh.
$$
This solution $u$ is actually $C^2$ by standard regularity results 
 for nonlinear uniformly elliptic equations \cite{KrylovE}.

To prove uniqueness of $\Hh$, suppose, by contradiction, that
$u_i$ and $\Hh_i$ ($i=1,2$) solve
$$
-\Delta u_i+H(D_xu_i,x)=\Hh_i.
$$
Suppose $u_1-u_2$ has a local maximum at $x_0$. Then
$D_xu_1=D_xu_2$ and $\Delta u_1\leq \Delta u_2$ at $x_0$.
Thus we conclude $\Hh_1\geq \Hh_2$. By symmetry $\Hh_1=\Hh_2$.

To prove that the viscosity solution is unique suppose, by contradiction, that
$u$ and $v$ are two distinct solutions (i.e., $u-v$ is non constant)
of
$$
-\Delta u+H(D_xu,x)=\Hh.
$$
We may assume
some small ball centered at the origin
of radius $\gamma$
does not contain
any maximizer of $x_0$ of $u(x)-v(x)$ (otherwise,
for convenience, we may shift the coordinates).
Fix $\epsilon,\lambda>0$
and assume that 
$$
u(x)-v(x)-\epsilon e^{-\lambda |x|^2}
$$
has a local maximum at $x_{\epsilon,\lambda}$. First observe that
$x_{\epsilon,\lambda}$ is uniformly bounded and 
by passing to a subsequence, if necessary, we may assume
$x_{\epsilon,\lambda}\rightarrow x_0$ as $\epsilon\rightarrow 0$.
At $x_{\epsilon,\lambda}$ we have
$$
D_xu=D_xv-2\epsilon\lambda x e^{-\lambda |x|^2}
$$
and
$$
\Delta u\leq \Delta v-2\epsilon\lambda e^{-\lambda |x|^2}+
\epsilon\lambda^2 |x|^2 e^{-\lambda |x|^2}).
$$
Since
$$
0=-\Delta (u-v)+H(D_xu,x)-H(D_xv,x),
$$
we have
$$
0\geq -2\epsilon\lambda e^{-\lambda |x_{\epsilon,\lambda}|^2}
+\epsilon\lambda^2 |x_{\epsilon,\lambda}|^2
e^{-\lambda |x_{\epsilon,\lambda}|^2}+O(\epsilon\lambda).
$$
Observe that for $\epsilon$ small enough
$|x_{\epsilon,\lambda}|>\frac{\gamma}{2}$.
Dividing by $\epsilon e^{-\lambda |x_{\epsilon,\lambda}|^2}$
and letting $\epsilon\rightarrow 0$
we observe that
$$
\frac{\gamma^2}{4}|\lambda|^2-C\lambda\leq 0.
$$
Therefore sending $\lambda\rightarrow\infty$ yields a contradiction.
\cqd

\begin{pro}
$\Hh$ can be estimated independently of $\sigma$ by
$$
\inf_x H(0,x)\leq \Hh\leq \sup_x H(0,x).
$$
\end{pro}
\Pf
Suppose $u$ has a minimum at $x_0$. Then $-\frac{\sigma^2}{2}\Delta u(x_0)
\leq 0$ and $D_xu(x_0)=0$. Thus
$$
\Hh=-\frac{\sigma^2}{2}\Delta u(x_0)+H(D_xu,x_0)\leq H(0,x_0)\leq
\sup_x H(0,x).
$$
The other estimate is similar.
\cqd

Finally we recall that standard estimates for controlled diffusions
\cite{FS} also yield that $u$ is semiconcave (with semiconcavity
constant independent of $\sigma$) and Lipschitz (also independently
of $\sigma$).

%Next we prove that a periodic semiconcave function is
%globally Lipschitz.
%\begin{lem}
%Suppose $u$ is periodic and semiconvex. Then $u$ is
%Lipschitz. 
%\end{lem}
%\Pf First we assume that $u$ is $C^2$. Let $x$ be an arbitrary point
%and let $p=D_xu(x)$. Choose $x_0$ to be a point for which
%$D_xu(x_0)=0$ and $y_k=x_k+k$ with $k\in\Zz^n$. Then
%$$
%(x-y_k)p\geq -C|x-y_k|^2\geq -C|x_0-y_k|^2-C
%$$
%Thus
%$$
%k p\geq -C|k|^2-C-C|p|,
%$$
%for all $k\in\Zz^n$. There exists $\delta>0$ such that
%for some $k_0\in\Zz^n$ we have
%$$
%k_0 p>\delta |p|.
%$$
%Thus for all $\lambda\in\Zz$
%$$
%\lambda k_0 p+\frac{k_0 p}{\delta}+ C\lambda^2+C\geq 0.
%$$
%Thus $(k_0 p)^2-4(C+\frac{k_0 p}{\delta})C\leq 0$ which implies
%that $k_0 p\leq C$ and therefore
%$$
%|p|\leq C,
%$$
%i.e. $u$ is Lipschitz.
%\cqd
%\begin{cor}
%$u$ is Lipschitz
%\end{cor}

\section{Equivalence between weak and strong problems}
\label{seceq}

The next task is to prove that the value $H^\star$, computed by
considering a infimum over measures in $\Mn_0$ is the same as
$$
\Hh=-\inf_{\mu\in \Mm_0}(h_1^*(\mu)-h_2^*(\mu)).
$$
A useful characterization of $\Hh$
 is:
\begin{teo}
$\Hh$
is the unique value for which the equation
\begin{equation}
\label{cell}
-\frac{\sigma^2}{2}\Delta u+H(D_xu,x)=\Hh
\end{equation}
has a periodic viscosity solution.
\end{teo}
\Pf We know from theorem \ref{solcel} that there is a single number $\Hh$
for which (\ref{cell}) admits a periodic viscosity solution $u$.
We can use that solution to build a Markov feedback strategy to control
the diffusion:
$$
dx=-D_pH(D_xu,x)dt+\sigma dw.
$$
To this diffusion it corresponds a measure $\mu\in\Mm_0\cap\Mm_1$
for which
$$
\int Ld\mu=-\Hh.
$$
Thus
$$
-\Hh\geq \inf_{\Mm_0\cap\Mm_1}\int Ld\mu.
$$

Conversely, let $u(x)$ be a solution of (\ref{cell}) and
assume that
$$-\Hh<\inf_{\Mm_0\cap\Mm_1}\int Ld\mu.$$
Then for any control strategy $u^*$
and corresponding process $x^*$
and all large enough T
$$
\frac{1}{T}E\int_0^T L(x^*,u^*)>-\Hh+\epsilon.
$$
Thus
$$
u(x)=\inf_{v}E\int_0^T L(x,v)+\Hh dt+u(x^*(T))
>\epsilon T+\min_x u(x),
$$
which is a contradiction for T sufficiently large since $u$ is bounded.
\cqd

\begin{teo}
$H^\star$ is the unique value for which the equation
$$
-\frac{\sigma^2}{2}\Delta u+H(D_xu,x)=H^\star
$$
has a periodic viscosity solution.
\end{teo}
\Pf First suppose $u$ is a periodic viscosity solution of
$$
-\frac{\sigma^2}{2}\Delta u+H(D_xu,x)=\Hh.
$$
Then we claim that there is no smooth function $\psi$ with
$$
-\frac{\sigma^2}{2}\Delta u+H(D_x\psi,x)<\Hh.
$$
Indeed, if this were false, we could choose a point $x_0$ at which
$u-\psi$ has a local minimum. At this point we would have
$$
-\frac{\sigma^2}{2}\Delta \psi+H(D_x\psi,x_0)\geq \Hh,
$$
by the viscosity property. Hence $H^\star\geq \Hh$, by
theorem \ref{nlm}.

To prove the other inequality consider a standard mollifier
$\eta_\epsilon$ and define $u_\epsilon=\eta_\epsilon*u$, 
in which $*$ denotes convolution.
Then
$$
-\frac{\sigma^2}{2}\Delta u_\epsilon+
H(D_xu_\epsilon,x)\leq \Hh+h(\epsilon,x),
$$
where
$$
h(\epsilon,x)=\sup_{|p|\leq R}\sup_{|x-y|\leq \epsilon}
|H(p,x)-H(p,y)|,
$$ 
where $R$ is a bound on the Lipschitz constant of $u$.
Let
$$H^\epsilon=\Hh+\sup_x h(\epsilon,x).$$
$u_\epsilon$ satisfies
$$
-\frac{\sigma^2}{2}\Delta u_\epsilon+
H(D_xu_\epsilon,x)\leq H^\epsilon.
$$
Thus $H^\star\leq \lim_{\epsilon\rightarrow 0} H^\epsilon=\Hh$.
Hence $H^\star=\Hh$.
\cqd

This proof holds even when $\sigma=0$, for $\sigma\neq 0$ since
$u$ is $C^2$, the mollification step is unecessary.

\begin{cor}
We have
$$
\inf_{\mu\in \Mn_0}(h_1^*(\mu)-h_2^*(\mu))=
\inf_{\mu\in \Mm_0}(h_1^*(\mu)-h_2^*(\mu)).
$$
\end{cor}
\Pf
Our previous results show that we can construct a probability
measure $\mu$ on $\Mm_0$
such that
$$
\int Ld\mu=-\Hh=\inf_{\mu\in \Mn_0}(h_1^*(\mu)-h_2^*(\mu)).
$$
Since $\Mm_0\subset \Mn_0$ this completes the proof.
\cqd

%Actually, we can prove that $\Mn_0$ is the weak-$*$
%closure of
%$\Mm_0$ (to a related proof see \cite{gdshgdhs}).
%To see this
%note that if $\Mm_0\neq \Mm_2$ then there would be a linear
%functional $L$ and a number $\alpha$ such that
%$$
%\int Ld\mu>\alpha\qquad \forall \mu\in\Mm_0
%$$
%and $\nu\in\Mm_2$ for which 
%$$
%\int Ld\nu>\alpha.
%$$

%Observe that if $u$ is a viscosity solution of the cell problem
%then
%$$
%A^v u+L\geq -\Hh
%$$
%at all points where the previous expression makes sense.
%Furthermore we can construct a sequence of smooth periodic
%$\psi_n$ such that
%$$
%A^v\psi+L\geq -\Hh-\delta_n
%$$
%Thus for any measure in $\Mm_2\cap\Mm_1$ we have
%$$
%\int Ld\mu\geq -\Hh
%$$

\section{Properties of Stochastic Mather measures}
\label{sec6}

In this section we
study general properties of Stochastic Mather measures.
First we prove that the stochastic Mather measure is supported
in the graph $(x,-D_pH(D_xu,x))$ for any $u$ viscosity solution of
(\ref{cell}). Then we show that the projection of this measure in
the $x$ axis has a density that satisfies an elliptic
partial differential equation.

\begin{teo}
Any stochastic Mather measure is supported in the graph $(x,-D_pH(D_xu,x))$
for any $u$ viscosity solution of
(\ref{cell}).
\end{teo}
\Pf
% Suppose $\mu$ is not
%supported on $(x,-D_pH(D_xu,x))$.
Recall that for any $v$ we have
$$
-\frac{\sigma^2}{2}\Delta u-vD_xu-L(x,v)\leq \Hh
$$
with strict inequality unless $v=-D_p\Hh(D_xu,x)$. Note that
$$
\int -\frac{\sigma^2}{2}\Delta u-vD_xu d\mu=0
$$
and
$$
-\int L(x,v)d\mu= \Hh,
$$
Thus $\mu$ is supported on $(x,-D_pH(D_xu,x))$,
otherwise we would have
$$
-\int L(x,v)d\mu< \Hh
$$
which would be a contradiction.
\cqd

Since any stochastic Mather measure is supported on a graph, a natural
question is whether its projection in the $x$ coordinates has a density.
The answer to this question is affirmative and we prove that
this density is the solution of an elliptic partial differential equation.

\begin{teo}
Let $\mu$ be a stochastic Mather measure.
Let $\nu$ denote the projection of $\mu$
in the $x$ coordinates. Then
$\nu=\theta(x)dx$ for some density
$\theta\in W^{1,2}$. Furthermore $\theta$ is a weak solution of
\begin{equation}
\label{measeq}
-\nabla(\theta v(x))+\frac{1}{2}\sigma^2\Delta\theta=0.
\end{equation}
for $v=-D_pH(D_xu,x)$.
\end{teo}
\Pf Recall that for any smooth and periodic $\phi(x)$ 
$$
\int \frac{\sigma^2}{2}\Delta\phi+v(x)D_x\phi d\nu=0.
$$
Let $\eta_\epsilon$ be a standard mollifier,
$\phi_\epsilon=\eta_\epsilon*\eta_\epsilon*\nu$
and $\nu_\epsilon=\eta_\epsilon*\nu$.
Note that $\nu_\epsilon$ is a bounded periodic $C^\infty$ functions
(the bounds may depend on $\epsilon$).
Then
$$
0=\int \frac{\sigma^2}{2}|\nabla \nu_\epsilon|^2dx-
\int v(x)D_x(\phi_\epsilon)d\nu.
$$
Thus
$$
\int \frac{\sigma^2}{2}|\nabla \nu_\epsilon|^2dx
=\int (D_x\nu_\epsilon) \eta_\epsilon*(v\nu)dx
$$
Note that
\begin{align}\notag
&|\int (D_x\nu_\epsilon) \eta_\epsilon*(v\nu)dx|\leq
|\int (D_x\nu_\epsilon) v(x)\nu_\epsilon dx|+\\\notag &\qquad
+|\int (D_x\nu_\epsilon) (\eta_\epsilon*(v\nu)-v(x)\nu_\epsilon )dx|.
\end{align}
The first term on the right-hand side can be estimated by
$$
\frac{\gamma}{2} \int |D_x\nu_\epsilon|^2 dx+\frac{C}{\gamma}
\int |\nu_\epsilon|^2dx,
$$
for any small $\gamma>0$.
To estimate the second term observe that since $v$ is
Lipschitz
$$
|\eta_\epsilon*(v\nu)-v(x)\nu_\epsilon|\leq
\int \eta_\epsilon(x-y)|v(x)-v(y)|d\nu(y)\leq
C\epsilon \nu_\epsilon.
$$
Thus
$$
|\int D_x\nu_\epsilon (\eta_\epsilon*(v\nu)-v(x)\nu_\epsilon )dx|
\leq \frac{\gamma}{2} \int |D_x\nu^\epsilon|^2 dx+\frac{C\epsilon}{\gamma}
\int |\nu_\epsilon|^2dx,
$$
therefore we conclude that
$$
\int |D_x\nu_\epsilon|^2dx\leq C\int |\nu_\epsilon|^2dx,
$$
uniformly in $\epsilon$.
Now observe that $\nu_\epsilon\geq 0$ and
$$
\int \nu_\epsilon dx=1.
$$
If $\int |\nu_\epsilon|^2dx$ were unbounded then we could normalize it
defining $\alpha_\epsilon =\gamma_\epsilon \nu_\epsilon$
with $\int |\alpha_\epsilon|^2dx=1$ and $\gamma_\epsilon\rightarrow 0$.
Since $\alpha_\epsilon \in W^{1,2}$ uniformly, through some
subsequence it converges in $L^2$ to some $\alpha\in L^2$ with
$\int |\alpha|^2dx=1$. However $\alpha\geq 0$ and $\int \alpha=0$
which is a contradiction.
Therefore we must
have $\nu_\epsilon\in W^{1,2}$ uniformly in $\epsilon$.
Thus through some subsequence $\nu_\epsilon\rightharpoonup \theta$
for some $\theta\in W^{1,2}$. Thus $d\nu=\theta(x)dx$.
Consequently, $\theta$ is a weak solution of
$$
-\nabla(\theta v(x))+\frac{1}{2}\sigma^2\Delta\theta=0.
$$
\cqd

Observe that equation (\ref{measeq}) is a non-symmetric zero eigenvalue
problem. It is well known \cite{Protter} that 
$$
-\nabla(\theta v(x))+\frac{1}{2}\sigma^2\Delta\theta=\lambda\theta
$$
has a principal eigenvalue $\lambda$ with positive eigenfunction $\theta$.
To see that $\lambda=0$ just observe that
$$
0=\int \frac{\sigma^2}{2}\Delta \theta-
\nabla\left(v(x)\theta(x)\right)=\lambda\int \theta.
$$
Since $\theta$ is non-negative we get $\lambda=0$.

The previous theorem yields several important identities that we will use in
the next section. First
define $\Hh(P)$ to be number for which
\begin{equation}
\label{celp}
-\frac{\sigma^2}{2}\Delta u+H(P+D_xu,x)=\Hh(P)
\end{equation}
has a periodic viscosity solution
$u(x,P)$ (note that
$u$ may not be continuous in $P$). The function $\Hh(P)$ is
convex in $P$ and so twice differentiable for almost
every $P$.

\begin{pro}
For any $\phi(x)$ periodic
\begin{equation}
\label{id1}
-\int D_x\phi D_pH\theta dx+\frac{\sigma^2}{2}
\int \Delta \phi \theta dx=0.
\end{equation}
Furthermore
\begin{equation}
\label{id2}
\int D_xH\theta dx=0.
\end{equation}
Finally, for any $P$ and $P'$,
\begin{equation}
\label{id3}
(P'-P)\int D_pH\theta dx \leq \Hh(P')-\Hh(P),
\end{equation}
in particular if $\Hh$ is differentiable
$$
\int D_pH\theta dx=D_P\Hh(P).
$$
\end{pro}
\Pf
Observe that (\ref{id1}) follows from
$$
0=\int \phi \left(\nabla(\theta D_pH)+\frac{1}{2}\sigma^2\Delta\theta\right)dx
$$
by integration by parts.

Let $\eta_\epsilon$ be a standard mollifier and let
$u_\epsilon=\eta_\epsilon*u$. Then
$$
-\frac{\sigma^2}{2}\Delta u_\epsilon+\eta_\epsilon* H(D_xu,x)=\Hh.
$$
Differentiate the previous identity with respect to $x_i$:
$$
-\frac{\sigma^2}{2}\Delta D_{x_i}u_\epsilon
+\eta_\epsilon*\left(H_{p_j}D_{x_jx_i}u+H_{x_i}\right)=0.
$$
Since $H_{p_j}$ is Lipschitz in $x$ we have
$$
\eta_\epsilon*(H_{p_j}D_{x_jx_i}u)=H_{p_j}D_{x_jx_i}u_\epsilon+O(\epsilon).
$$
Note also that
$$
\int \left(-\frac{\sigma^2}{2}\Delta D_{x_i}u_\epsilon+
H_{p_j}D_{x_jx_i}u_\epsilon\right)\theta dx=0
$$
since $D_{x_i}u_\epsilon$ is smooth and periodic.
Thus
$$
\int \eta_\epsilon*H_{x_i}\theta dx\rightarrow 0
$$
as $\epsilon\rightarrow 0$. Since $\eta_\epsilon*H_{x_i}\rightarrow H_{x_i}$
almost everywhere we conclude
$$
\int H_{x_i}\theta dx=0,
$$
which proves (\ref{id2}).

To prove the last part of the proposition, note that 
\begin{align}\notag &
H(P'+D_xu(x,P'),x)-H(P+D_xu(x,P),x)\geq\\\notag & \qquad
D_pH(P+D_xu(x,P),x) \left[P'+D_xu(x,P')-P
-D_xu(x,P)\right].
\end{align}
Let $w=u(x,P')-u(x,P)$. Note that
$$
\int \left(-\frac{\sigma^2}{2}\Delta w+
D_pH(P+D_xu(x,P),x)D_xw
\right)\theta dx=0.
$$
Thus
$$
(P'-P)\int  D_pH(P+D_xu(x,P),x)\leq \Hh(P')-\Hh(P),
$$
as required.
\cqd

\section{Regularity estimates}
\label{sec7}

In this section we prove $L^2$-type regularity estimates
for the solution of (\ref{celp}). These estimates are
expressed using the invariant measure. A major advantage
is that it is possible to prove $L^2(\theta)$ estimates for the difference
quotient $\left|D_xu(x+y)-D_xu(x)\right|$ that do not depend
on $\sigma$ explicitly whereas pointwise
or $L^2$ estimates with respect to Lebesgue measure depend on $\sigma$. 
Therefore our estimates extend up to the case $\sigma=0$,
for a careful study of this case consult
\cite{EGom1}, \cite{Gomes3}, and
\cite{Gomes2}.

\begin{teo}
Suppose $u$ solves (\ref{celp}) and $y\in\Rr^n$. Then
\begin{equation}
\label{est1}
\int \left|D_xu(x+y)-D_xu(x)\right|^2\theta dx\leq C|y|^2.
\end{equation}
Furthermore, if $\Hh(P)$ is twice differentiable at $P$ then
\begin{equation}
\label{est2}
\int \left|D_xu(x,P)-D_xu(x,P')\right|^2\theta dx\leq C|P-P'|^2,
\end{equation}
for $|P-P'|$ sufficiently small.
\end{teo}
\Pf Note that
$$
-\frac{\sigma^2}{2}\left[\Delta u(x+y)-\Delta u(x)\right]+H(D_xu(x+y),x+y)
-H(D_xu(x),x)=0.
$$
Since $H$ is convex,
$$
\begin{array}{l}
H(D_xu(x+y),x+y)-H(D_xu,x)\geq\gamma\left|D_xw\right|^2+\\
+D_pH(D_xu(x),x)D_xw
+D_xH(D_xu(x),x)y+O(y^2),
\end{array}
$$
with $w=u(x+y)-u(x)$.
Integrating with respect to $\theta dx$ to obtain
$$
\gamma\int \left|D_xw\right|^2\theta dx\leq
C|y|^2,
$$
since
$$
\int \left[D_pH(D_xu(x),x)D_xw-\frac{\sigma^2}{2}\Delta w
\right]\theta dx=0,
$$
and
$$
\int D_xH(D_xu(x),x)\theta dx=0.
$$

Similarly, let $w= u(x,P')-u(x,P)$ and assume $\Hh(P)$ is
twice differentiable at $P$. Then
$$
\begin{array}{l}
D_P\Hh(P)(P'-P)+C|P-P'|^2\geq\\
-\frac{\sigma^2}{2}\Delta w
+H(P'+D_xu(x,P'),x)-H(P+D_xu(x,P),x)
\end{array}
$$
Note that
\begin{align}
\notag
&\qquad H(P'+D_xu(x,P'),x)-H(P+D_xu(x,P),x)\geq\\\notag &\geq
D_pH(P+D_xu(x,P),x) (P'-P+D_xw)+\gamma
\left|P'-P+D_xw\right|^2.
\end{align}
Thus
\begin{equation}
\label{step}
\gamma \int \left|P'+D_xu(x,P')-P-D_xu(x,P)\right|^2\theta dx\leq
C|P-P'|^2,
\end{equation}
since
$$
(P'-P)\int D_pH(P+D_xu(x,P),x)\theta dx=(P'-P)D_P\Hh, 
$$
and
$$
\int -\frac{\sigma^2}{2}\Delta w+D_pH(P+D_xu(x,P),x) D_xw=0.
$$
From (\ref{step}) we have (\ref{est2}).\cqd

In the next theorem we prove that if $\Hh$ is strictly convex
in a neighborhood of a point $P$
then the map $(x,P)\rightarrow P+D_xu(x,P)$ is non-degenerate.
In the non-random case this result is extremely
important since it proves
the invariant sets $(x,P+D_xu)$ change with $P$, see \cite{Gomes3}
for a detailed discussion.

\begin{teo}
Suppose $\Hh$ is strictly convex at a neighborhood of a point $P$. Then
\begin{equation}
\label{est3}
\int \left|P+D_xu(x,P)-P'-D_xu(x,P')\right|^2\theta dx\geq C|P-P'|^2,
\end{equation}
for $|P-P'|$ sufficiently small.
\end{teo}
\Pf
let $w= u(x,P')-u(x,P)$ and assume $\Hh(P)$ is
strictly convex in a neighborhood of $P$. Then
\begin{align}\notag &
D_P\Hh(P)(P'-P)+C|P-P'|^2\leq\\\notag&\qquad\leq 
-\frac{\sigma^2}{2}\Delta w
+H(P'+D_xu(x,P'),x)-H(P+D_xu(x,P),x)
\end{align}
Note that
\begin{align}
\notag
&\qquad H(P'+D_xu(x,P'),x)-H(P+D_xu(x,P),x)\leq\\\notag &\leq
D_pH(P+D_xu(x,P),x) (P'-P+D_xw)+\Gamma
\left|P'-P+D_xw\right|^2.
\end{align}
Thus
$$
\Gamma \int \left|P'+D_xu(x,P')-P-D_xu(x,P)\right|^2\theta dx\geq
C|P-P'|^2
$$
since
$$
(P'-P)\int D_pH(P+D_xu(x,P),x)\theta dx=(P'-P)D_P\Hh 
$$
and
$$
\int -\frac{\sigma^2}{2}\Delta w+D_pH(P+D_xu(x,P),x) D_xw=0.
$$
\cqd

In the case $\sigma=0$ it is possible to prove $L^\infty$-estimates
on $D^2_{xx}u$ on the support of $\theta$ \cite{EGom1}.
However, this is not
the case for $\sigma> 0$, at least with estimates independent on
$\sigma$. Indeed, if $D^2_{xx}u$ were uniformly bounded in $\sigma$
then $u_\sigma$ would converge uniformly, through some
subsequence as $\sigma\rightarrow 0$,
to a function $u$, viscosity solution of
$$
H(D_xu,x)=\Hh.
$$
But then $u$ would be both semiconvex and semiconcave and we know that,
in general, $u$ is only semiconcave.
However, some regularity exists, as was remarked in
section \ref{ApA},
namely one-sided bounds on $D^2_{xx}u$ (semiconcavity)
that do not depend on $\sigma$.

\section{Explicit Formulas and Examples}
\label{sec8}

In this section we discuss several formulas for both $\Hh$ and
invariant measures.
%In section (\ref{seceq}) we derived a representation formula
%for $\Hh$ as
%$$
%\Hh=\inf_{\phi} \sup_x-\frac{\sigma^2}{2} \Delta\phi + H(D_x\phi,x).
%$$
%Define
%$$
%\mathcal{H}(\phi)= \sup_x-\frac{\sigma^2}{2} \Delta\phi + H(D_x\phi,x)
%$$
%then the function $\Hh$ is a convex function. Then
%$$
%\Hh=\inf_\phi \mathcal{H}(\phi)
%$$
%which is an $L^\infty$ calculus of variation problem.
%This representation formula is quite useful in computing $\Hh$ numerically.
%For instance in figure ... we can see the graph of $\Hh(P)$ for
%$$
%-\frac{1}{2}\Delta u+\frac{(P+D_xu)^2}{2}+\cos(x)=\Hh(P).
%$$ 
%
The next proposition shows that
given the solution $u(x,P)$ it is possible to
compute the density $\theta$ (under smoothness assumptions).
not of the invariant measure but of a time-reversed version.
\begin{pro}
Assume $u(x,P)$ is a smooth solution of (\ref{celp}). Then
$$
\theta=\det (I+D^2_{xP}u)
$$
is a solution of a time-reversed version of (\ref{measeq}):
\begin{equation}
\label{22a}
\frac{\sigma^2}{2}\Delta \theta+\nabla(\theta v(x))=0.
\end{equation}
\end{pro}
\Pf Let $v(x,P)=Px+u(x,P)$. Then
\begin{equation}
\label{v1}
-\frac{\sigma^2}{2}\Delta v+H(D_xv,x)=\Hh(P).
\end{equation}
The claim
is that $\theta=\det D^2_{xP}v$ solves
$$
-\frac{\sigma^2}{2}\Delta \theta+\nabla(\theta D_pH(D_xv,x))=0.
$$
Differentiate (\ref{v1}) with respect to $P_i$ to get
$$
-\frac{\sigma^2}{2}\Delta v_{P}+D_pHD^2_{xP}v=D_P\Hh
$$
Note that $D_pHD^2_{xP}v=(D^2_{xP}v)^T D_pH$ and multiply the 
previous identity by the cofactor matrix $\cof D_{xP}v$
$$
-\frac{\sigma^2}{2}\cof D_{xP}v\Delta v_{P}+\det D_{xP}v D_pH=
\cof D_{xP}vD_P\Hh.
$$
Observe that $\cof D_{xP}v$ is divergence free \cite{Evans1} and so
$$
\cof D_{xP}v\Delta v_{P}=\nabla (\det D_{xP}v).
$$
Therefore
$$
-\frac{\sigma^2}{2}\nabla (\det D_{xP}v)+\det D_{xP}v D_pH
=\cof D_{xP}vD_P\Hh.
$$
By applying $\nabla$ to the previous identity we have
$$
-\frac{\sigma^2}{2}\Delta \theta+\nabla(\theta D_pH)=0.
$$
\cqd
%
%Note that if $\sigma=0$ then $\det (I+D^2_{xP}u)$ is in fact an
%invariant measure \cite{EGom2}.

Now we turn our attention to
the special case
$$
H(p,x)=\frac{p^2}{2}+V(x),
$$
with $V$ periodic. For this special Hamiltonian we will
present an alternative representation formula for $\Hh(P)$
as well as exhibit a (non-periodic) invariant measure.
This will follow some ideas of (\cite{Holland}).

Suppose $u$ is a periodic viscosity solution of
$$
-\frac{\sigma^2}{2}\Delta u+H(P+D_xu,x)=\Hh(P).
$$
Define
$$
\phi=e^{-\frac{Px+u}{\sigma^2}}.
$$
Then $\phi$ solves
$$
\frac{\sigma^4}{2} \Delta \phi+V(x)\phi=\Hh(P)\phi.
$$
Thus $\Hh$ is an eigenvalue of the operator
$\frac{\sigma^4}{2} \Delta \phi+V(x)\phi$.
Consider the related operator
$$
L\psi=e^{\frac{Px}{\sigma^2}}\frac{\sigma^4}{2} \Delta
(e^{-\frac{Px}{\sigma^2}}\psi)+V(x)\psi
=\frac{\sigma^4}{2} \Delta \psi-\sigma^2 P D_x\psi+
(V(x)+\frac{|P|^2}{2})\psi
$$
Then $\Hh$ is also an eigenvalue of $L$ with periodic boundary conditions.
%Define
%$$
%J[\psi]=-\int_{\Tt^n} \psi L\psi=\frac{\sigma^4}{2}\int_{\Tt^n} |\nabla\psi|^2
%-(V(x)+\frac{\sigma^4|P|^2}{2})|\psi|^2,
%$$
%since $\int_{\Tt^n} \psi P D_x\psi=0$.
\begin{pro}
$\Hh$ is the principal eigenvalue of $L$.
\end{pro}
\Pf The operator $L$ has a principal eigenvalue $\lambda$
with positive and periodic eigenfunction $\varphi$. Let $u=-\log \varphi$.
Then $u$ is smooth, periodic and satisfies the Hamilton-Jacobi
equation
$$
-\frac{\sigma^2}{2}\Delta u+H(P+D_xu,x)=\lambda.
$$
By uniqueness of $\Hh$ we have $\lambda=\Hh$.
\cqd

Finally we exhibit an invariant measure for this system. Although this is
not a probability measure (unless $P=0$).
\begin{pro}
Let $\theta = e^{-2\frac{Px+u}{\sigma^2}}$. Then $\theta$ is an
invariant measure.
\end{pro}
\Pf It suffices to check that
$$
\frac{\sigma^2}{2}\Delta\theta+\nabla((P+D_xu)\theta)=0.
$$
\cqd

\section{Asymptotics}
\label{sec9}

In this section we study the asymptotic behavior of the controlled
process $x(t)$. First we will do 
some formal calculations motivated by the case $\sigma=0$
\cite{EGom1}, \cite{Gomes3}, \cite{Gomes2}.
Define
$$
X=x+D_Pu.
$$
Then
$$
dX=dx+D^2_{Px}u dx+\frac{\sigma^2}{2}D^3_{Pxx}udt
$$
Thus, since
$dx=-D_pHdt+\sigma dw$,
$$
dX=\left(-D_pH(I+D^2_{Px}u)+\frac{\sigma^2}{2}D^3_{Pxx}u\right)dt+
\sigma\left(1+D^2_{Px}u\right)dw.
$$
Note that $-D_pH(I+D^2_{Px}u)+\frac{\sigma^2}{2}D^3_{Pxx}u=-D_P\Hh$
and so
$$
E\left(X(t)-X(0)\right)=-D_P\Hh t.
$$
\begin{teo}
Suppose $\Hh$ is differentiable at $P$. Then
$$
\lim_{t\rightarrow\infty}E\frac{x(t)}{t}=-D_P\Hh.
$$
\end{teo}
\Pf
Let $u$ be a viscosity solution of (\ref{celp}).
Let $v^*$ be an optimal control such that
$$
u(x,P)=E\int_0^tL(x,v^*)+Pv^*+\Hh(P)+u(x(t),P).
$$
Then
$$
u(x,P')\leq E\int_0^tL(x,v^*)+P'v^*+\Hh(P')+u(x(t),P').
$$
Subtracting these two equations 
$$
C\leq E\int_0^t (P'-P)v^*+\Hh(P')-\Hh(P).
$$
Thus
$$
E\int_0^t v^*=-D_P\Hh t+O(1)
$$
since $dx=v^*dt+\sigma dw$ we have
$$
E\int_0^t v^*=E\int_0^tdx= Ex(t).
$$
\cqd

\section{Convergence as $\sigma\rightarrow 0$}
\label{sec10}

In this last section we prove that stochastic Mather measures
converge to a Mather measure as the diffusion rate $\sigma$
vanishes.

Let $\Hh_\sigma$ be the unique number
for which
\begin{equation}
\label{qw}
-\frac{\sigma^2}{2}\Delta u_\sigma+H(D_xu_\sigma,x)=\Hh_\sigma
\end{equation}
has a periodic viscosity solution $u_\sigma$. The bounds
on $\Hh_\sigma$ obtained in section \ref{ApA} imply that
through some subsequence $\Hh_\sigma\rightarrow \Hh$
as $\sigma\rightarrow 0$, for some number $\Hh$. Since $u_\sigma$
is uniformly Lipschitz in $\sigma$, through some subsequence
$u_\sigma\rightarrow u$ uniformly. Standard stability results on
viscosity solutions imply that $u$ is a viscosity solution of
$$
H(D_xu,x)=\Hh.
$$

Let $\mu_\sigma$ be a stochastic Mather measure 
associated with (\ref{qw}). Since the support of $\mu_\sigma$
is bounded independently of $\sigma$ we can extract
a weakly convergence subsequence $\mu_\sigma\rightharpoonup \mu$
and $\int d\mu=1$.
Note that
$$
-\Hh_\sigma=\int Ld\mu_\sigma\rightarrow \int Ld\mu=-\Hh
$$
Furthermore, for any smooth function $\phi(x)$
$$
0=\int \frac{\sigma^2}{2}\Delta \phi+vD_x\phi d\mu_\sigma
\rightarrow \int vD_x\phi d\mu.
$$
Thus $\mu$ satisfies
$$
\int Ld\mu =-\Hh
$$
with the constraints $\int d\mu=1$, and $\int vD_x\phi d\mu=0$.
Thus $\mu$ is a Mather measure.

\bibliographystyle{alpha}

\bibliography{duality2}

\end{document}